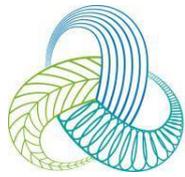 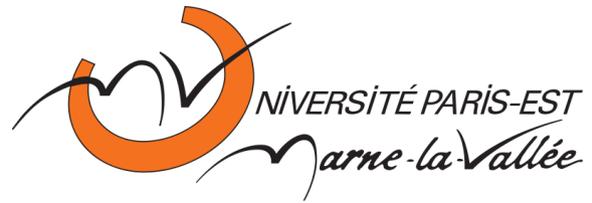

# Multi-Anticipative Piecewise Linear Car-Following Model

Nadir Farhi[*], Habib Haj-Salem[†] & Jean-Patrick Lebacque[‡]

Université Paris-Est, IFSTTAR, GRETTIA, F-93166 Noisy-le-Grand, France.

[*] nadir.farhi@ifsttar.fr (Corresponding author).
[†] habib.haj-salem@ifsttar.fr
[‡] jean-patrick.lebacque@ifsttar.fr

**Abstract**

We propose in this article an extension of the piecewise linear car-following model to multi-anticipative driving. As in the one-car-anticipative model, the stability and the stationary regimes are characterized thanks to a variational formulation of the car-dynamics. We study the homogeneous driving case. We show that in term of the stationary regime, the multi-anticipative model guarantees the same macroscopic behavior as for the one-car-anticipative one. Nevertheless, in the transient traffic, the variance in car-velocities and accelerations is mitigated by the multi-anticipative driving, and the car-trajectories are smoothed. A parameter identification of the model is made basing on NGSIM data and using a piecewise linear regression approach.

# INTRODUCTION

We present a multi-anticipative car-following traffic model, where drivers control their velocities by taking into account the positions and the velocities of many cars ahead. In basic car-following models [3, 6, 7, 10], the car dynamics are described by stimulus- response equations that express the control process of drivers. Each driver accelerates or decelerates depending on his speed, and on the relative speed and the inter-vehicular distance with respect to the driver of the car ahead.

Multi-anticipative car-following models are often extensions of existing one car anticipative models. Since 1968, Bexelius [2] extended the model of Chandler et al. [3] to the multi-anticipative case. Lenz et al. [12] extended the model of Bando et al. [1]. More recently, Hoogendoorn et al. [11] have extended the model of Helly [8, 9].

The model we present here is an extension of the piecewise linear car following model [4, 5]. It is a first-order discrete-time model where the car-velocities are given in function of the inter-vehicular distances. We show here that the variational formulation made in [4, 5] holds also for the multi-anticipative extension we propose in this article. That is to say that the car-dynamics are again interpreted as dynamic programming equations associated to stochastic optimal control problems of Markov chains, as in [4, 5]. Thanks to that formulation, we are able to characterize the stability of the car-dynamics and to calculate the stationary regimes.

In the transient traffic, some qualitative results are obtained from our model. First, we simulated the car dynamics on a one lane-road of about 10,000 meters, where we imposed the trajectories of a given number of leaders, and simulated the trajectories of all the followers, basing on our model, and varying the number of leaders taken into account in anticipation. The effect of our anticipation modeling on the transient traffic is thus shown. As expected, we observed that, as the number of leaders taken into account in anticipation, increases, the car trajectories are smoothed, the distant followers are slowed down in unstable traffic phases, but they retrieve their non-anticipative trajectories once the traffic is stabilized.

Second, we proposed a parameter identification approach for our multi-anticipative model, and showed, in a basic example, how this method is applied. We rely on NGSIM data of vehicle trajectories on a segment of U.S. highway 101. The purpose of this part is rather to present the process of parameter identification, since the data used here are not exhaustive to draw conclusions, and the model is likely to be improved. However, some interesting observations are made on the data sample considered. For example, we observed that the scatter plot for space headway and average velocities seems to be easy to approximate by a piecewise linear curve in the case where the space headway is computed with respect to two leaders, compared to the case where the space headway is computed with respect to only one leader.

In the remainder of this introduction, we give a short review on the car-following modeling, we fix notations, and we introduce the next sections. We use the notations $t$ for time (discrete or continuous), $x$ for distance (car positions), and $n$ for the number of cars. The cars are numbered such that the first car (car number 1) is the leader. We consider the variables

- $x(n, t)$: the cumulative traveled distance of car $n$, from time zero to time $t$.
- $y(n, t)$: the inter-vehicular distance $x(n - 1, t) - x(n, t)$.
- $v(n, t)$: the velocity of car $n$ at time $t$.

Car-following models are often based on a behavioral law $V_e$ (equilibrium speed spacing function), that gives, at the equilibrium traffic, the velocity $v$ of a car $n$ as a function of the inter-vehicular distance $y$ between cars $n$ and $n - 1$. It is then assumed that the law $V_e$ also holds on the transient traffic. A kind of general form of first-order car-following models can then be derived as done for the macroscopic first-order modeling (LWR models [13,14]).

$$\dot{v}(n,t) = V'_e(y(n,t))\Delta v(n,t), \tag{1}$$

where $\Delta v(n,t) = v(n-1,t) - v(n,t)$.

The simplest form for the equilibrium speed-spacing function $V_e(y)$ is the linear one $V_e(y) = \alpha y + \beta$, where $\alpha$ and $\beta$ are parameters. In this is case, the linear car-following model [3, 10] is obtained:

$$v(n,t) = \alpha y(n,t) + \beta, \tag{2}$$

The model (2) is, of course, not satisfactory. Another well known car-following model is due to Gazis, Herman, and Rothery [6].

$$\dot{v}(n,t+T) = a\, v(n,t+T)^p \frac{v(n-1,t) - v(n,t)}{(x(n-1,t) - x(n,t))^l} \tag{3}$$

where a reaction time $T$ is considered, and where $p$ and $l$ are parameters. For example, if $p = 1$ and $l = 2$, the model can simply be obtained by taking $V_e(y) = b\exp{-a/y}$ in (1), satisfying $V'_e(y) = a\, V_e(y)/y^2$.

Other car-following models that do not necessarily match the form (1) exist, such as the models of Bando et al. [1] (optimal velocity model), Helly [8, 9], Treiber et al. [15] (intelligent driver model), etc. Bando et al. [1] proposed the optimal velocity model

$$\dot{v}(n,t+T) = \lambda(V_e(y(n,t)) - v(n,t+T)) \tag{4}$$

This model has been much studied recently for being easily analyzed with mathematical tools.

Helly [8, 9] considered the linear model

$$\dot{v}(n,t+T) = \alpha\, \Delta v(n,t) + \beta(y(n,t) - S_n) \tag{5}$$

where $\alpha$ and $\beta$ are parameters, and where $S_n$ is the desired distance which can be linear $S_n = S_0 + Tv_n$, with $S_0$ the minimum gross distance between two cars.

We base here on the piecewise linear car-following model proposed in [5], where the behavioral law $V_e$ is approximated with a (min-max)-piecewise linear curve

$$V_e(y) = \min_{u \in U} \max_{w \in W} \{\alpha_{uw} y + \beta_{uw}\} \tag{6}$$

and where a one-car-anticipative discrete-time car-dynamics has been obtained

$$x(n,t+1) = x(n,t) + \min_{u \in U} \max_{w \in W} \{\alpha_{uw}(x(n-1,t) - x(n,t)) + \beta_{uw}\} \tag{7}$$

where $\alpha_{uw}$ and $\beta_{uw}$, for $(u,w) \in U \times W$, are parameters, and $U$ and $W$ are two finite sets of indices. The system (7) is also written, for the traffic of $\nu$ cars $1, 2, \ldots, \nu$ as follows.

$$x_n(t+1) = \min_{u \in U} \max_{w \in W} \{[M^{uw} x(t)]_n + c_n^{uw}\}, \quad 1 \leq n \leq \nu \tag{8}$$

where $M^{uw}$ and $c^{uw}$, for $(u,w) \in U \times W$, are matrices and column-vectors respectively.

Two cases have been distinguished in [5].

- The $\nu$ cars move on a ring road. In this case, $M^{uw}$ and $c^{uw}$ are given by

$$M^{uw} = \begin{pmatrix} 1-\alpha_{uw} & 0 & \cdots & \alpha_{uw} \\ \alpha_{uw} & 1-\alpha_{uw} & & 0 \\ \vdots & \ddots & \ddots & \vdots \\ 0 & 0 & \alpha_{uw} & 1-\alpha_{uw} \end{pmatrix}$$

and

$$c^{uw} = {}^t\left(\frac{\alpha_{uw}\nu}{d} + \beta_{uw}, \beta_{uw}, \ldots, \beta_{uw}\right)$$

The dynamics (8) is stable under the condition $\alpha_{uw} \in [0,1]$, and the behavior law is realized at the stationary regime

$$\bar{v} = \min_{u \in U} \max_{w \in W} \{\alpha_{uw}\bar{y} + \beta_{uw}\} \qquad (9)$$

where $\bar{v}$ denotes the asymptotic car-velocity (the same for all cars), and $\bar{y}$ denotes the average inter-vehicular distance in the ring road ($\bar{y} = 1/d$).

- The $\nu$ cars move on an "open" road, where the velocity $v_1(t)$ of the first car (the leader one) varies over time but is stationary. In this case, $M^{uw}$ and $c^{uw}$ are given by

$$M^{uw} = \begin{pmatrix} 1 & 0 & \cdots & 0 \\ \alpha_{uw} & 1-\alpha_{uw} & & 0 \\ \vdots & \ddots & \ddots & \vdots \\ 0 & 0 & \alpha_{uw} & 1-\alpha_{uw} \end{pmatrix}$$

and

$$c^{uw}(t) = {}^t(v_1(t), \beta_{uw}, \ldots, \beta_{uw})$$

Again, the dynamics (8) is stable under the condition $\alpha_{uw} \in [0,1]$ and the inverse behavior law at the stationary regime is obtained as follows.

$$\bar{y} = \max_{u \in U} \min_{w \in W} \frac{v_1 - \beta_{uw}}{\alpha_{uw}} \qquad (10)$$

where $v_1$ denotes the asymptotic velocity of the first car.

The dynamics (8) have been interpreted in [5], under the assumption $\alpha_{uw} \in [0,1], \forall (u,w) \in U \times W$, as a dynamic programming equation associated to a stochastic game on a controlled Markov chain; see [5] for more details.

## 1 ANTICIPATION MODELLING

We present in this section an extension of the model (7) to multi-anticipative traffic, where each car chooses its velocity depending on the inter-vehicular distance with respect to a given number m of cars ahead the considered car (multi-leaders). In order to situate our model with respect to the existing multi-anticipative models, and to explain the extension we do, let us first give a short review on multi-anticipative car-following models.

A straightforward multi-leader extension of the model of Chandler et al. [3] is the Bexelius model [2]

$$\dot{v}(n, t + T) = \sum_{j=1}^{m} \alpha_j \Delta v^{(j)}(n, t) \qquad (11)$$

where $\dot{v}$ denotes the acceleration, $\alpha_j, j = 1,2, \ldots, m$ are sensitivity parameters with respect to the $j$th car ahead, and where $\Delta v^{(j)}(n, t) = v(n - j, t) - v(n, t)$. This model is very simple but permits some mathematical analysis.

Hoogendoorn et al. [11] have noted the non convenience of the additive form of Bexelius model (11), and proposed the modification

$$\dot{v}(n, t + T) = \min_{1 \leq j \leq m} \alpha_j \Delta v^{(j)}(n, t) \qquad (12)$$

Hoogendoorn et al. [11] have also proposed a multi-anticipative generalization for the Helly model (5)

$$\dot{v}(n, t + T) = \sum_{j=1}^{m_1} \alpha_j \Delta v^{(j)}(n, t) + \sum_{j=1}^{m_2} \beta_j [\Delta x^{(j)}(n, t) - S^j(n)] \qquad (13)$$

where $\Delta x^{(j)}(n, t) = x(n - j, t) - x(n, t)$. Lenz et al. [12] have generalized the Bando model (4) as follows.

$$\dot{v}(n, t) = \sum_{j=1}^{m} \kappa_j \left\{ V_e \left( \frac{\Delta x^{(j)}(n, t)}{j} \right) - v(n, t) \right\} \qquad (14)$$

where $\kappa_j$ expresses the sensitivity with respect to the $j$th leader.

We propose here a multi-leader extension for the piecewise linear car-following model (7). We use a minimum form as in (12) (rather than an additive form as in (11)). Moreover, we use a uniform form for the sensitivity with respect to the inter-vehicular distance as in (14) (the inter-vehicular distance with respect to the $j$th leader is divided by $j$). We consider the dynamics:

$$x_n(t + 1) = x_n(t) + \min_{1 \leq j \leq m} (1 + \lambda)^{j-1} \min_{u \in U} \max_{w \in W} \left\{ \alpha_{uw} \left( \frac{x_{n-j}(t) - x_n(t)}{j} \right) + \beta_{uw} \right\} \qquad (15)$$

where $m$ is the number of leaders taken into account in anticipation, and $\lambda \geq 0$ is a discount parameter with respect to the leader index. The dynamics (15) can be written simply

$$x_n(t + 1) = x_n(t) + \min_{1 \leq j \leq m} \min_{u \in U} \max_{w \in W} \left\{ \alpha_{juw} \left( \frac{x_{n-j}(t) - x_n(t)}{j} \right) + \beta_{juw} \right\} \qquad (16)$$

where $\forall (u, w) \in U \times W, (\alpha_{juw})_j, 1 \leq j \leq m$ are increasing non negative sequences, and $(\beta_{juw})_j, 1 \leq j \leq m$ are increasing sequences.

The interpretation of the minimum operator with respect to the $j$th leader in (16) is that a car $n$ maximizes its velocity under the constraints

$$x_n(t + 1) - x_n(t) \leq \min_{u \in U} \max_{w \in W} \left\{ \alpha_{juw} \left( \frac{x_{n-j}(t) - x_n(t)}{j} \right) + \beta_{juw} \right\}, \quad 1 \leq j \leq m.$$

One consequence of anticipation in driving is that the information that a car $i$, for $i = n - 1, n - 2, \ldots, \max(1, n - m)$, decelerates at time $t$, is immediately transmitted to the car $n$ that reacts at time $t + 1$, instead of $t + n - i$. The discounting with respect to the leader indices, made by introducing the multiplicative term $(1 + \lambda)^{j-1}$, permits to favor closer leaders over distant ones. If $\lambda = 0$, then the cars respond equally to the stimulus of all the leaders $j$, with $j = 1, 2, \ldots, m$.

In the following section (section 2), we study the stability of the car dynamics (15), and characterize the existence of stationary regimes. Two cases are distinguished: traffic on a ring road, and traffic on an open road. In both cases, we give the asymptotic car positions when stationary regimes exist. The transient traffic for the car dynamics (15) is treated in section 3.

## 2 STABILITY ANALYSIS AND STATIONARY REGIMES

As in [4, 5], we consider $\nu$ cars moving on a 1-lane road without passing. We first study the case where the cars move on a ring road, and then explore the "open" road case.

### 2.1 Traffic on a ring road

The cars being moving on a ring road, the indices $n - j$, in the dynamics (15), are cyclic in the set $\{1, 2, \ldots, \nu\}$. The idea here is that the two minimum operators in (15) can be summarized in only one minimum operator, and then retrieve the one car-anticipative form for the dynamics. Let us denote by $Z$ the set of all pairs of indices $(j, u)$, with $1 \leq j \leq m$ and $u \in U$
$$Z = \{z = (j, u), 1 \leq j \leq m, u \in U\}.$$
The dynamics (15) is then written

$$x_n(t+1) = \min_{z \in Z} \max_{w \in W} \{[M^{zw} x(t)]_n + c_n^{zw}\}, \quad 1 \leq n \leq \nu \tag{17}$$

where the matrices $M^{zw} = M^{juw}$ and the column vectors $c^{zw} = c^{juw}$ are given as follows.

$$M^{juw} = \begin{pmatrix} 1 - \alpha_{juw}/j & 0 & \cdots & \alpha_{juw}/j & 0 & 0 \\ 0 & 1 - \alpha_{juw}/j & \ddots & \ddots & \alpha_{juw}/j & 0 \\ \vdots & \ddots & \ddots & \ddots & \ddots & \alpha_{juw}/j \\ \alpha_{juw}/j & \ddots & \ddots & \ddots & \ddots & \vdots \\ 0 & \alpha_{juw}/j & \ddots & \ddots & \ddots & 0 \\ 0 & 0 & \alpha_{juw}/j & \cdots & 0 & 1 - \alpha_{juw}/j \end{pmatrix}$$

and

$$c^{juw} = {}^t\!\left(\frac{(\alpha_{juw}/j)\nu}{d} + \beta_{juw}, \beta_{juw}, \ldots, \beta_{juw}\right).$$

The dynamics (18) have the same form as (8). It is then interpreted as a dynamic programming equation associated to a stochastic game on a controlled Markov chain. The stability is guaranteed under the condition $\alpha_{juw} \in [0,1], \forall (j, u, w) \in \{1, 2, \ldots, m\} \times U \times W$; see [4,5] for more details. The stationary regime is characterized by the additive eigenvalue problem

$$\bar{v} + x_n = \min_{1 \leq j \leq m} \min_{u \in U} \max_{w \in W} \{[M^{juw} x]_n + c_n^{juw}\}, \quad 1 \leq n \leq \nu \tag{18}$$

where $\bar{v}$ is the asymptotic car-velocity, the same for all cars, the vector $x$ is the asymptotic car-positions, given up to an additive constant. The following result gives a solution for the system (19).

**Theorem 1.** If $\forall (j, u, w) \in \{1, 2, \ldots, m\} \times U \times W, \alpha_{juw} \in [0,1]$, then the system (19) admits a solution $(\bar{v}, x)$ given by:
$$\bar{v} = \min_{u \in U} \max_{w \in W} \{\alpha_{1uw}\bar{y} + \beta_{1uw}\},$$
$$x = {}^t((v-1)\bar{y}, (v-2)\bar{y}, \ldots, \bar{y}, 0).$$

where $\bar{y} = 1/d$ is the average inter-vehicular distance in the ring road.

**Proof.** Following the same approach as in [4, 5], we obtain the following solution for the system (19).
$$\bar{v} = \min_{z \in Z} \max_{w \in W} \{\alpha_{zw}\bar{y} + \beta_{zw}\} = \min_{1 \leq j \leq m} \min_{u \in U} \max_{w \in W} \{\alpha_{juw}\bar{y} + \beta_{juw}\},$$
$$x = {}^t((v-1)\bar{y}, (v-2)\bar{y}, \ldots, \bar{y}, 0).$$

Then since $(\alpha_{juw})_j$ and $(\beta_{juw})_j$ are increasing sequences (with respect to $j$) and $\bar{y} \geq 0$, we have

$$\min_{1 \leq j \leq m} \min_{u \in U} \max_{w \in W} \{\alpha_{juw}\bar{y} + \beta_{juw}\} = \min_{u \in U} \max_{w \in W} \{\alpha_{1uw}\bar{y} + \beta_{1uw}\}.$$
∎

A particular case is important here, which is the non discounting case where $\lambda = 0$ in (15). In this case, the matrices $M^{juw}$ and the vectors $c^{juw}$ still depend on $j$, whilst the parameters $\alpha_{juw}$ and $\beta_{juw}$ are independent of $j$ for all $j \in \{1, 2, \ldots, m\}$. Thus the average car speed $\bar{v}$ coincides with the average car-speed obtained in the one car anticipative model (7):

$$\bar{v} = \min_{u \in U} \max_{w \in W} \{\alpha_{uw}\bar{y} + \beta_{uw}\},$$

where $\alpha_{uw} = \alpha_{juw}$ and $\beta_{uw} = \beta_{juw}$, $\forall j \in \{1, 2, \ldots, m\}$.

## 2.2 Traffic on an open road

We suppose here given the speed $v_1(t)$ of the first car over time, since the cars move on an open road. In addition, in order to analyze the stability and the stationary regime of the car dynamics, we assume that the velocity of the first car approaches a constant value $v_1$. That is to say that $\lim_{t \to +\infty} v_1(t) = v_1$. Moreover, the number of anticipation cars for the first $m$ cars cannot be $m$ (it is less than $m$). More precisely, the number of anticipation cars for a car numbered $n$ is $\min(n-1, m)$.

The dynamics (15) is written here

$$x_n(t+1) = \min_{1 \leq j \leq m} \min_{u \in U} \max_{w \in W} \{[M^{juw} x(t)]_n + c_n^{juw}\}, \quad 1 \leq n \leq v \tag{19}$$

where the matrices $M^{juw}$ and the vectors $c^{juw}$ are given by

$$M^{juw} : \begin{pmatrix} 1 & 0 & 0 & \cdots & & & & \cdots & 0 \\ 0 & 1 & 0 & 0 & \cdots & & & & 0 \\ \vdots & & \ddots & & & & & & \vdots \\ 0 & 0 & \cdots & 1 & 0 & 0 & \cdots & & 0 \\ \alpha_{juw}/j & 0 & \cdots & 0 & 1-\alpha_{juw}/j & \cdots & & & 0 \\ 0 & \alpha_{juw}/j & 0 & \cdots & 0 & 1-\alpha_{juw}/j & & & \vdots \\ \vdots & 0 & \ddots & 0 & \cdots & & \ddots & 0 & \\ 0 & \cdots & 0 & \alpha_{juw}/j & 0 & \cdots & & 0 & 1-\alpha_{juw}/j \end{pmatrix} \begin{matrix} 1 \\ 2 \\ \vdots \\ j \\ j+1 \\ j+2 \\ \vdots \\ v \end{matrix}$$

$$c^{juw} = \begin{pmatrix} v_1(t) \\ +\infty \\ \vdots \\ +\infty \\ \beta_{juw} \\ \beta_{juw} \\ \vdots \\ \beta_{juw} \end{pmatrix} \begin{matrix} 1 \\ 2 \\ \vdots \\ j \\ j+1 \\ j+2 \\ \vdots \\ v \end{matrix}$$

The entries $(M^{juw})_{ik}$ and $(c^{juw})_i$ for $i,k \leq j$ do not play any role in the car dynamics since $(c^{juw})_i = +\infty, \forall i \leq j$. Those entries correspond to anticipation of a car $i$ with respect to its $j^{th}$ leader that does not exist since $i \leq j$.

We assume that the velocity $v_1(t)$ of the first car reaches a fixed value $v_1$ at the stationary regime. The stationary regime is thus characterized as follows.

$$\bar{v} + x_n = \min_{1 \leq j \leq m} \min_{u \in U} \max_{w \in W} \{[M^{juw} x]_n + c_n^{juw}\}, \quad 1 \leq n \leq v \qquad (20)$$

The following result gives a solution for the system (21).

**Theorem 2.** For all $y \in R$ satisfying $\min_{u \in U} \max_{w \in W} (\alpha_{1uw} y + \beta_{1uw}) = v_1$, the couple $(\bar{v}, x)$ is a solution for the system (21), where $\bar{v} = v_1$ and $x$ is given up to an additive constant by

$$x = {}^t((v-1)y, (v-2)y, \ldots, y, 0). \qquad (21)$$

**Proof.** The proof is similar to that of Theorem 3 of [5]. Let $y \in R$ satisfying

$$\min_{u \in U} \max_{w \in W} (\alpha_{1uw} y + \beta_{1uw}) = v_1.$$

Let $x$ be given by (22). Then $\forall n \in \{1, 2, \ldots, v\}$ we have

$$\begin{aligned} \min_{1 \leq j \leq m} \min_{u \in U} \max_{w \in W} [M^{juw} x]_n + c_n^{juw} &= \min_{1 \leq j \leq m} \min_{u \in U} \max_{w \in W} (\alpha_{juw} y + \beta_{juw}) + x_n \\ &= \min_{u \in U} \max_{w \in W} (\alpha_{1uw} y + \beta_{1uw}) + x_n \\ &= v_1 + x_n. \end{aligned}$$

Moreover, the optimal strategy at the stationary regime is $(j, u, w) = (1, \bar{u}, \bar{w})$ such that $\alpha_{1\bar{u}\bar{w}} y + \beta_{1\bar{u}\bar{w}} = v_1$. Indeed

$$\begin{aligned} [M^{1\bar{u}\bar{w}} x]_n + c_n^{1\bar{u}\bar{w}} &= (\alpha_{1\bar{u}\bar{w}} y + \beta_{1\bar{u}\bar{w}}) + x_n \\ &= \min_{u \in U} \max_{w \in W} (\alpha_{1uw} y + \beta_{1uw}) + x_n \\ &= \min_{1 \leq j \leq m} \min_{u \in U} \max_{w \in W} (\alpha_{juw} y + \beta_{juw}) + x_n \\ &= \min_{1 \leq j \leq m} \min_{u \in U} \max_{w \in W} [M^{juw} x]_n + c_n^{juw}. \end{aligned}$$

∎

Theorem 2 gives the car-velocity at the stationary regime with one stationary configuration of cars (uniform configuration) and gives the optimal strategy for drivers at that regime. An important remark here is that the car-velocity obtained is the same as the one obtained for the one car anticipative model (7). Moreover, the "optimal strategy" of driving at the stationary regime in the case of multi-anticipative model is to drive by taking into account only one leader: $(\bar{j}, \bar{u}, \bar{w}) = (1, \bar{u}, \bar{w})$. That is, at the stationary regime, once the traffic is stabilized, it is not necessary for drivers to take into consideration more than one leader.

Moreover, the stationary configurations of cars in the two cases of one car anticipation and multi-anticipation models may coincide. This case is interpreted as follows. Even though the cars reduce their approach in the multi-anticipative dynamics (due to the minimum operator over the leader indices), comparing to their movement under the one car anticipative dynamics; as long as the cars approach the stationary regime, where the traffic is stabilized, they retrieve what they have lost on the transient regime. Therefore, by introducing the minimum on the multi-anticipative dynamics, the traffic becomes smoother, without decreasing the stationary car speed.

## 3 TRANSIENT TRAFFIC

We take here the same example presented in [5], which we adapt to multi-anticipative case, in order to make a comparison. We simulate the car-dynamics (15). We take as the time unit half a second (1/2 s), and as the distance unit 1 meter (m). The parameters of the model are the same as those of Example 1 of [5] (the parameters have been determined by approximating a given behavior law). More precisely, The behavior law considered here is approximated by the following piecewise linear curve of six segments.

$$\tilde{V}(y) = \max\{\alpha_1 y + \beta_1, \min\{\alpha_2 y + \beta_2, \alpha_3 y + \beta_3, \alpha_4 y + \beta_4, \alpha_5 y + \beta_5, \alpha_6 y + \beta_6\}\},$$

where the parameters $\alpha_i$ and $\beta_i$ for $i = 1, 2, \ldots, 6$ are given by

**TABLE 1 Approximation of the behavior law $\tilde{V}(y)$ with a piecewise-linear curve**

| Segments | 1 | 2 | 3 | 4 | 5 | 6 |
|---|---|---|---|---|---|---|
| $\alpha_i$ | 0 | 0.54 | 0.32 | 0.13 | 0.34 | 0 |
| $\beta_i$ | 0 | −8.1 | −1.47 | 6.11 | 10.6 | 14 |

We simulate the car dynamics on a one lane road of about 10,000 meter. We vary the number of leaders and the velocity of the first car (over time) in order to show the effect of multi-anticipation on the transient traffic. We give the results on Table 1. We see easily on that table that the trajectories are smoothed by anticipation. We notice here that the trajectory of the first car is the same for all views of Table 2. Although the number of leaders taken into account by drivers cannot exceed 5 in practice, we simulated here the car dynamics with anticipation with up to 100 leaders. We did it for curiosity, but it can be interesting in the case, for example, where one likes to study the traffic of communicating cars or automatic ones, etc.

## 4 PARAMETER IDENTIFICATION

We propose here a parameter identification method based on (piecewise) linear regression. Given measured data of the car positions and velocities on a given section, the method permits to determine the optimal parameters that match our model with the measured data. Since the dynamics of the model is simply the car-velocities given as functions of the inter-vehicular distances, then the optimal parameters would be the ones that approximate the scatter plot of instantaneous inter-vehicular distances and velocities.

**TABLE 2** Traffic on a 1-lane road. On the $x$-axis: time. On the $y$-axis: car-position. The number of cars taken into account in anticipation are 1, 5, 10, 20, 50 and 100. The length of the road is 10,000 meter. The total simulation time is 500 seconds

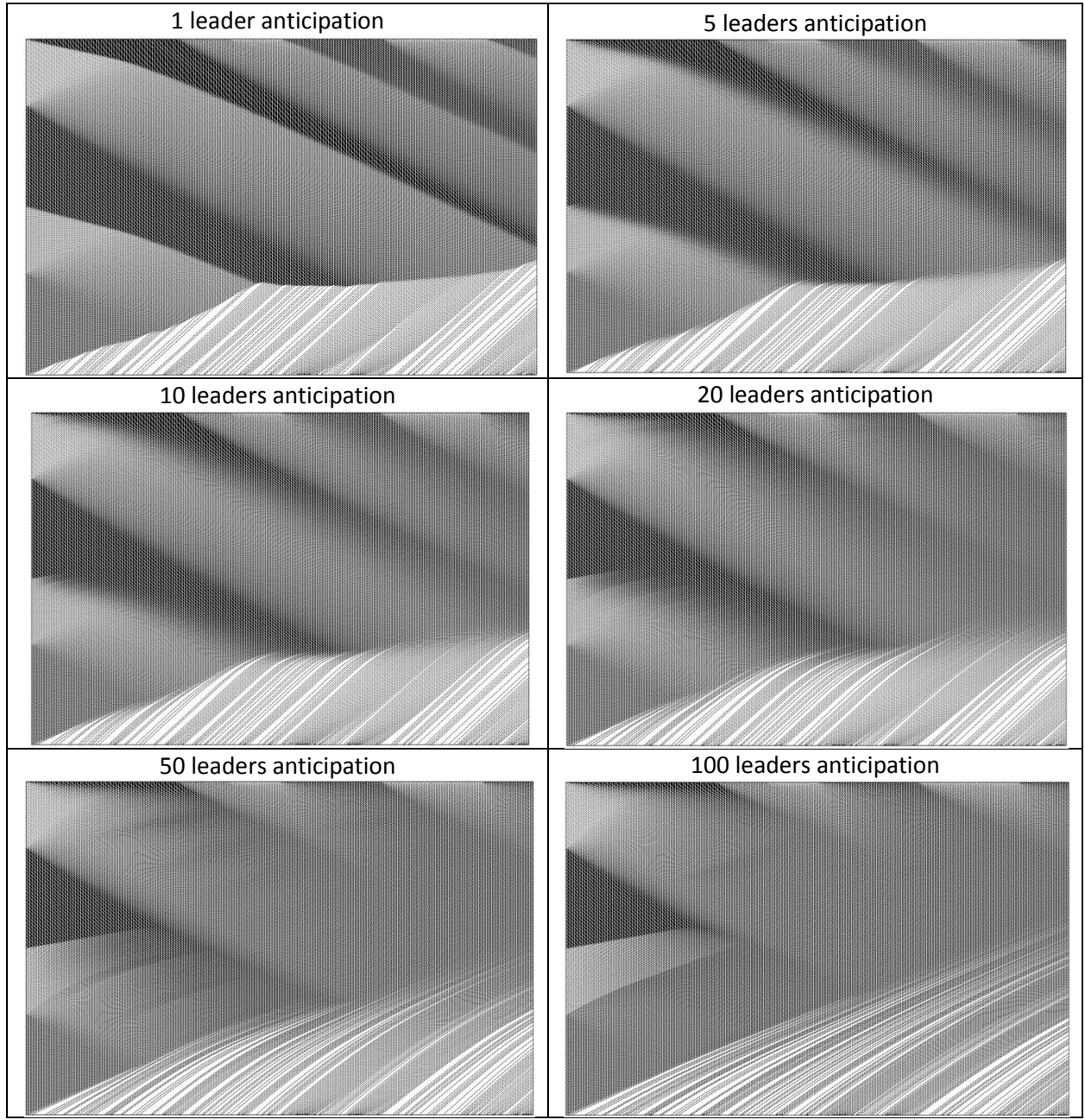

We denote by $\tilde{y}_{(m,\lambda)}$ the variable

$$\tilde{y}_{(m,\lambda)}(n,t) = \min_{1 \leq j \leq m} \frac{(1+\lambda)^{j-1}(x(n-j,t) - x(n,t))}{j}.$$

For fixed values of $m$ and $\lambda$ and with measured velocities $v(n,t)$ and inter-vehicular distances $x(n-j,t) - x(n,t)$ for every car $n$, we approximate the scatter plot $V(\tilde{y})$ by a piecewise-linear curve

$$V(\tilde{y}) = \min_{u \in U} \max_{w \in W} (\alpha_{uw} \tilde{y} + \beta_{uw}). \tag{22}$$

The min-max piecewise linear approximation (23) is based on a piecewise linear regression approach where the number of segments as well as the points where they intersect are determined optimally by deterministic dynamic programming. The optimization with respect to the parameters $m$ and $\lambda$ is done numerically by varying the two parameters in convenient intervals, and then determine the optimal parameters.

For fixed values of parameters $m$ and $\lambda$, the scatter plot $V(\tilde{y}_{(m,\lambda)})$ is approximated by linear regression on separated intervals. The intervals are determined by a dynamic programming approach. More precisely, we divide the axis of inter-vehicular distances $\tilde{y}$ into unity intervals $(\tilde{y}_i, \tilde{y}_{i+1})$. We start from the first interval, make a linear regression on that interval, and then decide for the second interval whether we make only one linear regression for the two intervals, or approximate the scatter plot on the second interval with another linear segment. Then we do the same for the third unity interval, and so on. We solve here the optimal control problem associated to the decision process. Note that the inter-vehicular distance plays the role of time in this decision process.

The decision at the unity interval $y_i$, denoted here by $r(y_i)$, is in $\{0,1\}$:

- $r(y_i) = 0$, if we decide to make one linear regression for the interval $y_i$ together with the intervals before it.

- $r(y_i) = 1$ if we decide to start a new linear regression from the unity interval $y_i$.

We define costs $k(y, N, r)$ in order to minimize regression errors and penalize large segmentation (limit the number of segments used in the approximation).

- $k(y_i, N, r)$ the error of regression at the unit interval $(y_i, y_{i+1})$, when the interval $(y_{i-N}, y_i)$ is wholy approximated by one segment, and when the decision $r$ is taken at the stage $y_i$.

Then the costs $k(y, N, r)$ are given by
- $k(y, N, r)$ = Linear regression error in $(y_{i-N}, y_{i+1})$ if $r = 0$,
- $k(y, N, r)$ = Linear regression error in $(y_i, y_{i+1}) + \varphi$ if $r = 1$,

where $\varphi$ is a penalization of starting a new linear regression. Then the following optimal control problem is solved

$$\min_{\gamma \in \Gamma} \sum_{y=0}^{y_{max}} k(y, N, r), \qquad (23)$$

where $\Gamma$ is the set of boolean strategies on $\{0,1, \dots, y_{max}\}$. That is $\gamma : \{0,1, \dots, y_{max}\} \ni y \mapsto r \in \{0,1\}$. The value function associated to (24) is

$$G(y, N) = \sum_{z=y}^{y_{max}} k(z, N, r) \qquad (24)$$

$G$ satisfies the dynamic programming equation

$$\begin{aligned} G(y_{max}) &= 0 \\ G(y_i, N) &= \min\{k(y_i, N, 0) + G(y_{i+1}, N+1), k(y_i, N, 1) + G(y_{i+1}, 1)\}. \end{aligned} \quad (25)$$

The parameter identification method is thus summarized as follows. For each couple of parameters $(m, \lambda)$, we calculate $\tilde{y}_{(m,\lambda)}$. Then the curves $V(\tilde{y}_{(m,\lambda)})$ are approximated by using the piecewise linear regression approach explained above (that is by solving the dynamic programming equation (26)). We obtain a total regression error for each approximation $(m, \lambda)$. Finally, we determine the best couple $(m, \lambda)$ that gives the minimal total error.

We show below, a first application of the identification method proposed above. We base here on NGSIM data of vehicle trajectories on a segment of U.S. Highway 101 (Hollywood Freeway) in Los Angeles, California. The data are collected between 7:50 a.m. and 8:05 a.m. on June 15, 2005. A preliminary analysis of the trajectories showed that, according to the multi-anticipative model presented here, the $n$th leaders' positions for $n > 3$ are redundant in the data considered here, even with a null discount parameter ($\lambda = 0$). That is to say that

$$\forall \lambda \geq 0, \forall n, t, \quad \tilde{y}_{(m,\lambda)}(n, t) \geq \tilde{y}_{(3,\lambda)}(n, t).$$

By consequent, we do not consider more than 3 leaders here ($m \in \{1,2,3\}$). Indeed, we can see from the dynamics (15), that, for $\lambda = 0$, taking into account more than one leader in anticipation does not change anything in the case where the traffic is accelerating, because the spacing to the $j$th leader is bigger than $j$ times the spacing to the first leader. That is

$$x_{n-j}(t) - x_n(t) \geq j \cdot (x_{n-1}(t) - x_n(t)), \quad j = 1, 2, \ldots$$

Nevertheless, this is realistic because, in case of accelerating traffic, anticipation is not significant.

First, we show in Table 3 the scatter plot for $\tilde{y}_{(m,0)}$ with the average car-velocities (the average over all cars), for $m \in \{1,2,3\}$. It seems that $V(\tilde{y}_{(2,0)})$ and $V(\tilde{y}_{(3,0)})$ (where we denote also by $V$ the average car-velocity over all cars) can easily be approximated with piecewise linear curves, comparing to $V(\tilde{y}_{(1,0)})$. We do not have justification for that presently.

**TABLE 3** Average car-velocity function of the inter-vehicular distance. $\lambda = 0$

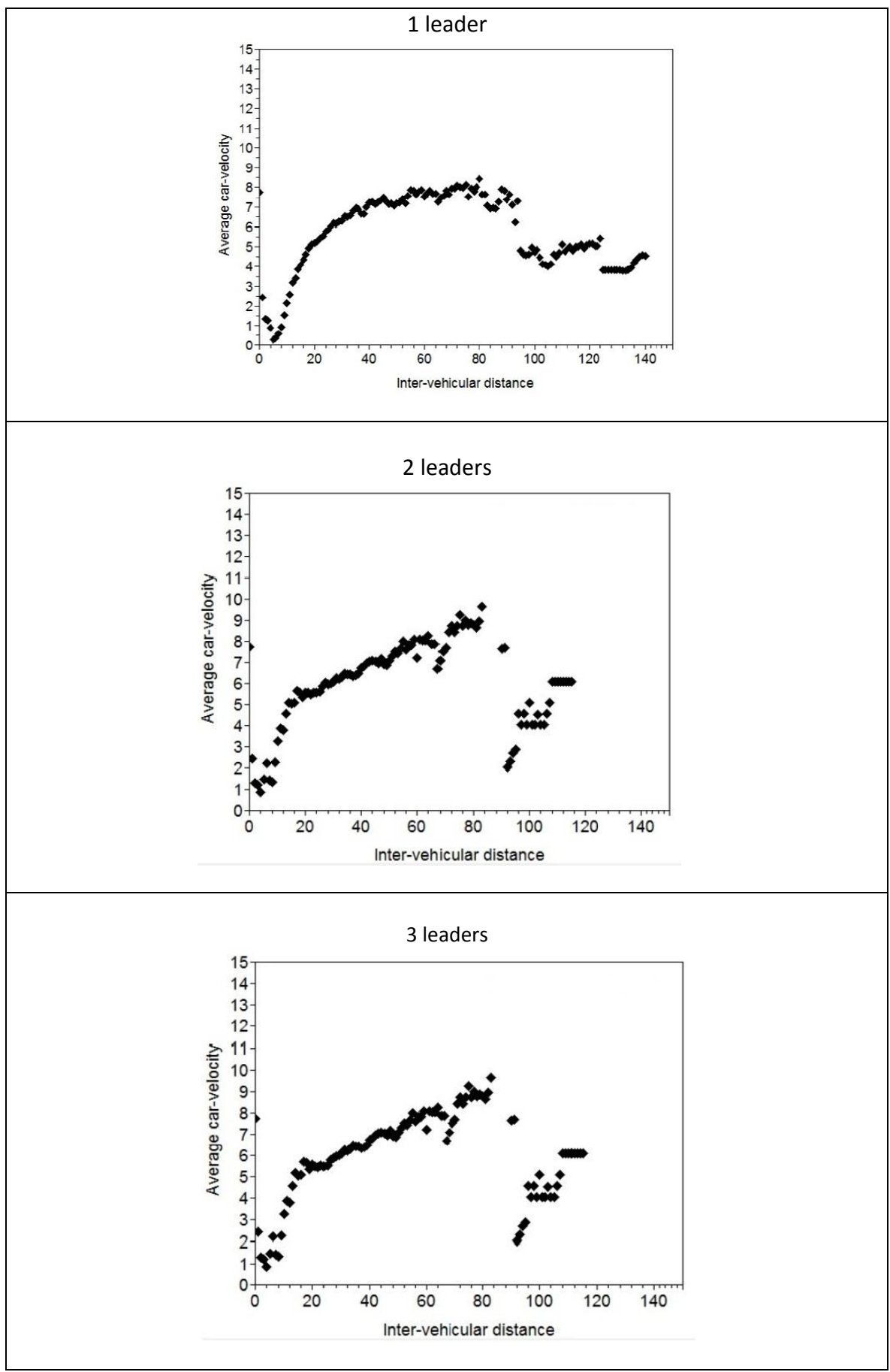

In Figure 1, we give the results of the parameter identification. The figure gives the total errors from the piecewise linear regression, obtained for varied values of parameters $m \in \{1,2,3\}$ and $\lambda \in [0, 5]$, increased by 0.1. The optimal parameters obtained here are $m = 2$ and $\lambda = 1.5$. As mentioned above, it is not worthy to consider values of $m$ that exceed 3 (for the data considered here). Also, for large values of $\lambda$, we retrieve the same results as if we take $m = 1$. Therefore, the optimization we did here is significant.

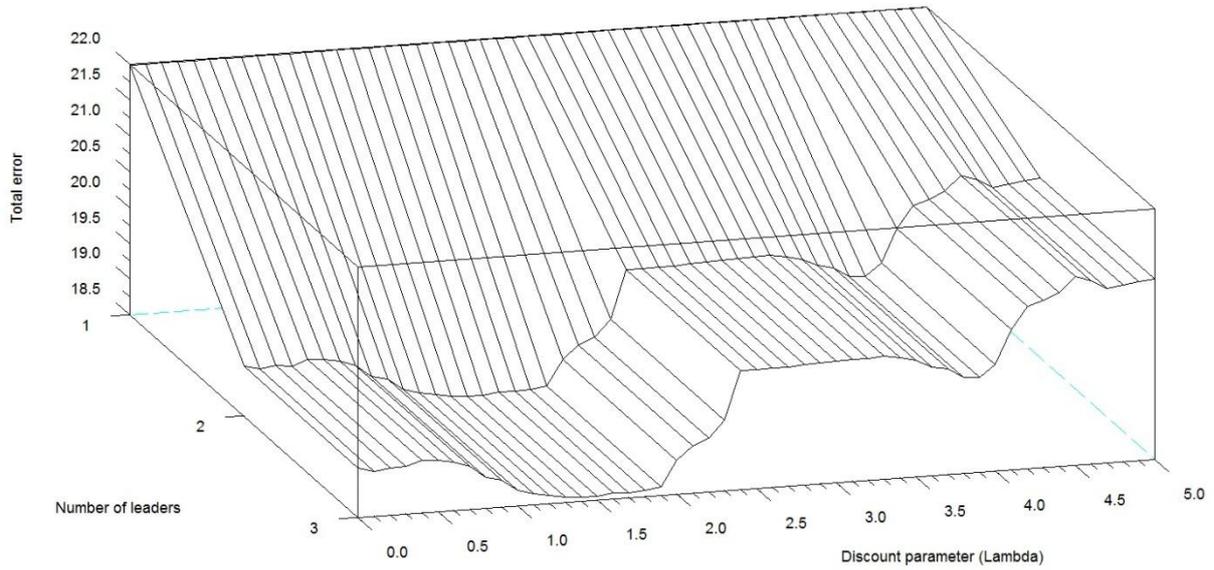

**FIGURE 1  Identification of parameters $m$ and $\lambda$**

The scatter plot for $\tilde{y}_{(2,1.5)}$ with $V$, is approximated by the following curve of four segments.

$$V(\tilde{y}) = \max\{0, \min\{0.38\tilde{y} - 1.90, 0.11\tilde{y} + 2.95, 10\}\}. \tag{26}$$

The approximation is shown in Figure 2.

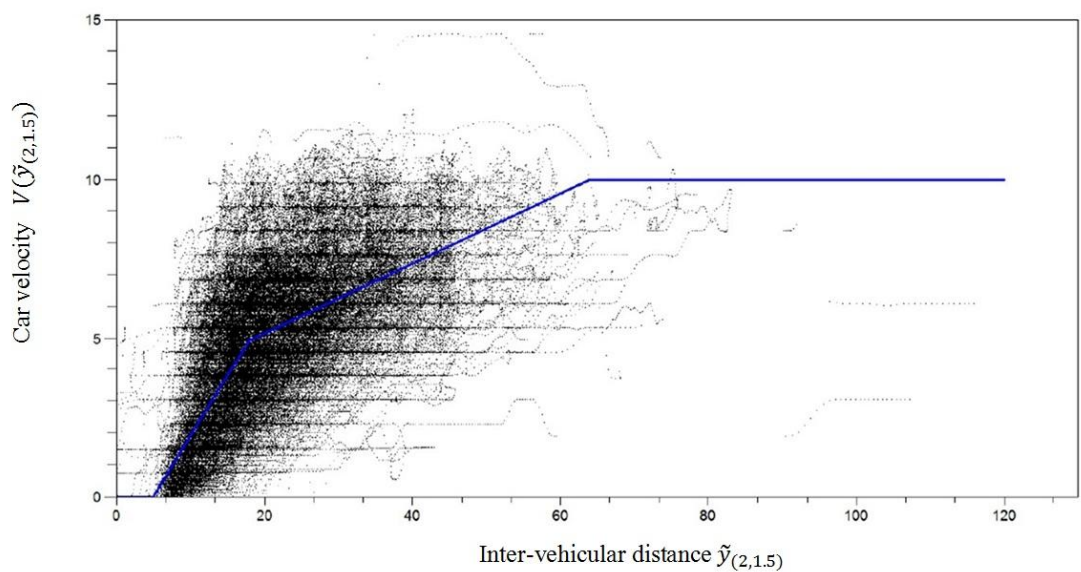

**FIGURE 2  Approximation of the law $V(\tilde{y}_{(2,1.5)})$ with a min-max piecewise linear curve. The number of leaders taken into account is $m = 2$. The discounting parameter is $\lambda = 1.5$**

As shown in Figure 2, it is not easy to identify one behavior for a big number of drivers. In fact the parameter identification should be made for each driver. The identification of the behavior law for a randomly selected driver has given the following result. The optimal parameters $(m, \lambda)$ are $m = 2$ (anticipation with two leaders) and $\lambda = 0$ (no discounting). The curve $V(\tilde{y}_{(2,0)})$ is approximated as follows.

$$V(\tilde{y}_{(2,0)}) \approx \max(0, \min(0.33\tilde{y}_{(2,0)} - 1.71, 8)).$$

The approximation is shown in Table 4, where we have shown also the approximation of the curve $V(\tilde{y}_{(1,0)})$.

$$V(\tilde{y}_{(1,0)}) \approx \max(0, \min(0.26\tilde{y}_{(1,0)} - 0.9, 8)).$$

As obtained by the parameter identification approach, we can see on Table 4 that the curve $V(\tilde{y}_{(2,0)})$ is well fitted by a piecewise-linear curve, comparing to the curve $V(\tilde{y}_{(1,0)})$. Note that the analytical results on the stability of the dynamics and the stationary regimes, presented in section 2 do not necessarily hold for the case of several driver behaviors. We will consider in a next step, the driver heterogeneity.

**TABLE 4   Parameter identification of one driver behavior**

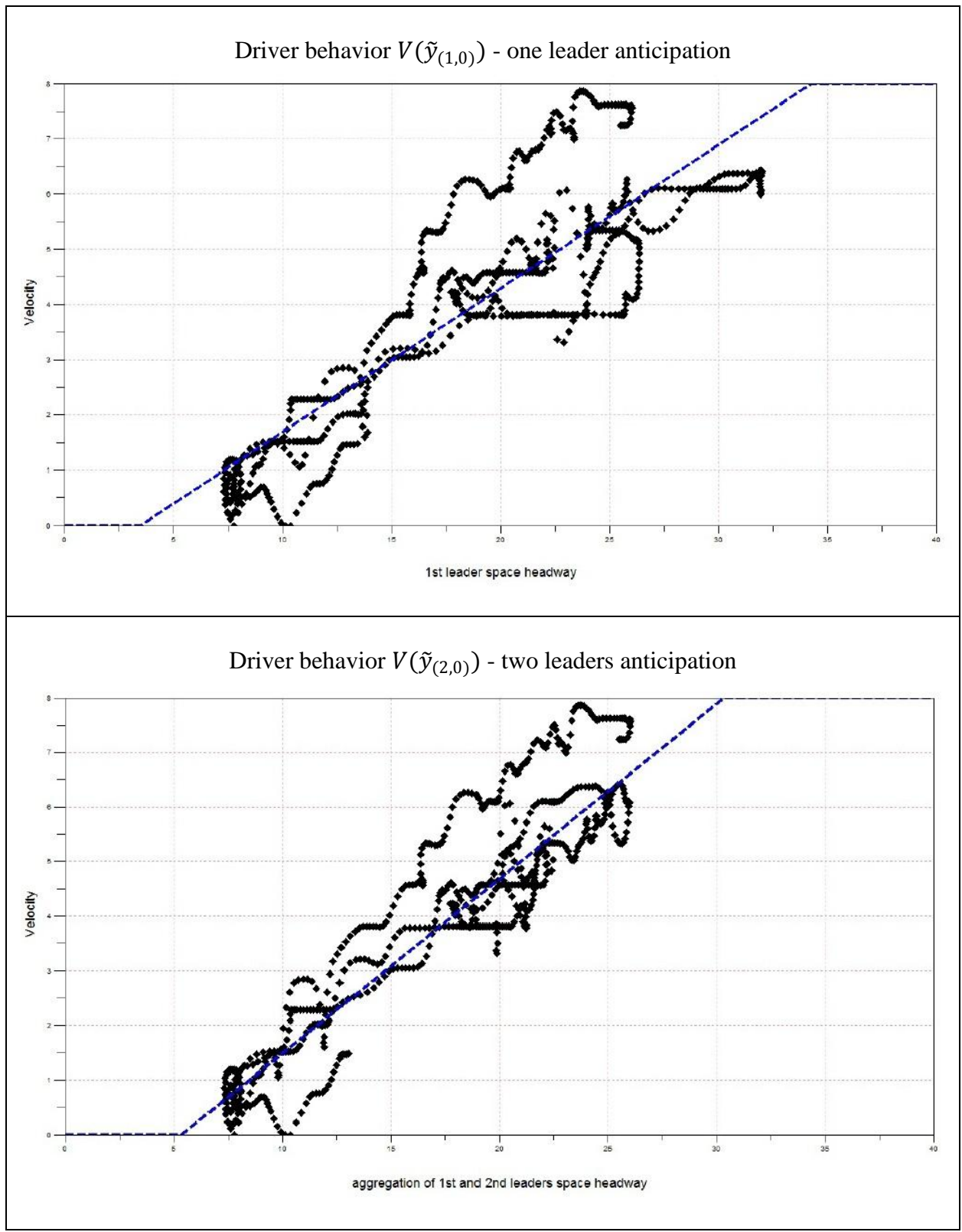

**5 CONCLUSION**

We presented in this article an extension of the piecewise linear car following model to multi-anticipative driving. The minimum form used for taking into account more than one leader, with the discounting parameter used to favor closest leaders over distant ones, seem to be convenient to recapture the main characteristics of the anticipative traffic. We have, in particular, shown that the trajectories are smoothed by the anticipation on the transient traffic, without affecting the stationary regimes. That is, the anticipation do not slow down the stationary traffic. The minimum form used for anticipation allows us to make the same variational formulation as done in the one-car anticipative traffic, and by this, it allows us to characterize the stability of the car-dynamics and calculate the stationary regimes. The identification test we made here is only a first step of analyzing the proposed model. In-depth analyses on exhaustive data should be done in the future, in order to improve the modeling approach presented in this article. In particular, we shall try to extend the model in a way that it takes into consideration heterogeneity in driving.

**REFERENCES**

[1] M. Bando, K. Hasebe, A. Nakayama, A. Shibata, and Y. Sugiyama. Dynamical model of traffic congestion and numerical simulation. Physical Review E, 51(2), 1995.

[2] S. Bexelius. An extended model for car-following. Transportation Research, 2(1):13– 21, 1968.

[3] R. E. Chandler, R. Herman, and E. W. Montroll. Traffic dynamics: Studies in car following. Operations Research, 6:165–184, 1958.

[4] Nadir Farhi. Modelisation Minplus et Commande du Trafic de Villes Regulieres. PhD thesis, University Paris 1 Panthon-Sorbonne, 2008.

[5] Nadir Farhi. Piecewise linear car-following modeling. arXiv:1107.5869, 2011.

[6] D C Gazis, R Herman, and R W Rothery. Nonlinear follow-the-leader models of traffic flow. Operations Research, 9(4):545–567, 1961.

[7] Denos C. Gazis, Robert Herman, and Renfrey B. Potts. Car-Following Theory of Steady-State Traffic Flow. Operations Research, 7(4):499–505, 1959.

[8] W. Helly. Simulation of bottlenecks in single lane traffic flow. In International Symposium on the Theory of Traffic Flow, 1959.

[9] W. Helly. Simulation of bottlenecks in single lane traffic flow. in theory of traffic flow. Elsevier Publishing Co., pages 207–238, 1961.

[10] R Herman, E W Montroll, R B Potts, and R W Rothery. Traffic dynamics: Analysis of stability in car following. Operations Research, 7(1):86–106, 1959.

[11] Serge P Hoogendoorn, Saskia Ossen, and M Schreuder. Empirics of multianticipative car-following behavior. Transportation Research Record, 1965:112120, 2006.

[12] H. Lenz, C. K. Wagner, and R. Sollacher. Multi-anticipative car-following model. European Physical Journal B, 7:331–335, 1999.

[13] J. Lighthill and J. B. Whitham. On kinematic waves ii: A theory of traffic flow on long, crowded roads. Proc. Royal Society, A229:281–345, 1955.

[14] P. I. Richards. Shock waves on the highway. Operations Research, 4:42–51, 1956.

[15] Martin Treiber, Ansgar Hennecke, and Dirk Helbing. Congested traffic states in empirical observations and microscopic simulations. PHYSICAL REVIEW E, 62:1805, 2000.